\documentclass[11pt]{amsart}%
\usepackage{amssymb,latexsym}
\usepackage{graphicx}
\usepackage{hyperref}
\usepackage{amsmath}
\usepackage{amsfonts}
\usepackage{amssymb}
\usepackage{float}%
\pdfoutput=1
\newtheorem{theorem}{Theorem}[section]
\newtheorem{proposition}[theorem]{Proposition}

\newtheorem{preremark}{Remark}
\newenvironment{remark}    {\begin{preremark}\rm}{\end{preremark}}
\newtheorem{lemma}[theorem]{Lemma}
\newtheorem{preexample}{Example}
\newenvironment{example}    {\begin{preexample}\rm}{\end{preexample}}

\newtheorem{prop}[theorem]{Proposition}

\theoremstyle{definition}
\theoremstyle{remark}

\renewcommand{\to}{\rightarrow}

\begin{document}
\title[Stellar subdivisions and Stanley-Reisner rings]{Stellar subdivisions and Stanley-Reisner rings of Gorenstein complexes}
\author{Janko B\"ohm}
\address{Janko B\"ohm\\
Department of Mathematics\\
University of Kaiserslautern\\
Erwin-Schr\"odinger-Str.\\
67663 Kaiserslautern\\
Germany}
\email{boehm@mathematik.uni-kl.de}
\author[Stavros~A.~Papadakis]{Stavros~Argyrios~Papadakis}
\address{Stavros~Argyrios~Papadakis, Centro de An\'{a}lise Matem\'{a}tica, Geometria e Sistemas Din\^{a}micos,
Departamento de Matem\'atica, Instituto Superior T\'ecnico, Av. Rovisco Pais,
1049-001 Lisboa, Portugal}
\email{papadak@math.ist.utl.pt}
\thanks{J. B. supported by DFG (German Research Foundation) through Grant BO3330/1-1.
S.~P. is a participant of the Project PTDC/MAT/099275/2008, a member of CAMGSD
(IST/UTL), and was supported from the Portuguese Funda\c{c}\~ao para a
Ci\^encia e a Tecno\-lo\-gia (FCT) under research grant SFRH/BPD/22846/2005 of POCI2010/FEDER}

\begin{abstract}
Unprojection theory analyzes and constructs complicated commutative
rings in terms of simpler ones. Our main result is that, on the algebraic
level of Stanley--Reisner rings, stellar subdivisions of 
Gorenstein* simplicial complexes correspond to unprojections of type
Kustin--Miller. As an application of our methods we study the 
minimal resolution of Stanley--Reisner rings 
associated to stacked polytopes, recovering results of 
Terai, Hibi, Herzog and Li Marzi. 
\end{abstract}
\maketitle

\section{Introduction}

\label{intro}

Stanley--Reisner rings of simplicial complexes form an important class of
commutative rings whose theory has provided spectacular applications to
combinatorics; see \cite{Sta} and \cite[Chapter 5]{BH} \cite{MS}. The
Stanley--Reisner ring of a simplicial complex $\Delta$, defined as the
quotient of a polynomial ring by a certain ideal, depends only on the
combinatorics of $\Delta$. Given a combinatorial operation on $\Delta$ which
produces another simplicial complex, it is natural to ask how the
Stanley--Reisner ring of the new complex is related to that of $\Delta$.
Stellar subdivision, which is one of the simplest ways to subdivide a
simplicial complex, is such an operation. It has been used successfully, for
instance, to give a method for transforming the boundary of a polytope into
that of any other polytope of the same dimension by operations which preserve
interesting invariants \cite{ES}, to construct polytopes whose $f$-vectors, or
flag $f$-vectors, span a certain `Euler' or `Dehn-Sommerville' space
\cite[Chapter 9]{Gr} \cite{BB} and to construct simplicial polytopes with
prescribed face lattices \cite{Lee, Si}.

On a different tone, unprojection theory aims to analyze and construct
commutative rings in terms of simpler ones. The first kind of unprojection
which appeared in the literature is that of type Kustin--Miller, studied
originally by Kustin and Miller \cite{KM} and later by Reid and the second
author \cite{P1, PR}. Starting from a codimension one ideal $I$ of a
Gorenstein ring $R$ such that the quotient $R/I$ is Gorenstein, Kustin--Miller
unprojection uses the information contained in $\operatorname{Hom}_{R}(I,R)$
to construct a new Gorenstein ring $S$ which is `birational' to $R$ and
corresponds to the `contraction' of $V(I)\subset\operatorname{Spec}R$. It has
been used in the classification of Tor algebras in Gorenstein codimension $4$
\cite{KM2}; in the birational geometry of Fano $3$-folds \cite{CM, CPR}; in
the study of Mori flips \cite{BrR}; in the study of algebraic surfaces of
general type \cite{NP}, \cite{NPi}; in the construction of weighted K3
surfaces and Fano $3$-folds \cite{Al}, \cite{BrKeRe}; and in the construction
of Calabi--Yau $3$-folds of high codimension \cite{B,NP2}. A general
discussion of unprojection theory and its applications is contained in
\cite{R1}, while a precise general definition of unprojection is proposed in
\cite{P2}. The Kustin--Miller unprojection and the associated complex
construction has been implemented in the package \textsc{KustinMiller}
\cite{BP2} for the computer algebra system \textsc{Macaulay2} \cite{GS}.

The main objective of this paper is to show that the Stanley--Reisner rings of
stellar subdivisions of a Gorenstein* simplicial complex $\Delta$
can be constructed from the Stanley--Reisner ring of $\Delta$ by unprojections
of type Kustin--Miller. As an application, we inductively calculate the
minimal graded free resolution of the Stanley--Reisner rings of the boundary
simplicial complexes of stacked polytopes, recovering results by 
Terai and Hibi \cite{TH2} and Herzog and Li Marzi \cite{HM}.

To state our main result, we need to introduce some notation and terminology
(see Section \ref{sec:pre} for more details). We denote by $k[ \Delta]$ the
Stanley--Reisner ring of a simplicial complex $\Delta$ with coefficients in a
fixed field $k$. Recall that $\Delta$ is said to be Gorenstein* over $k$ if
$k[ \Delta]$ is Gorenstein and given a vertex $i$ of $\Delta$ there 
exists $\sigma \in \Delta$ such that  $\sigma \cup \{ i \}$ is
not a face of $\Delta$.

Given a face $\sigma$ of $\Delta$, we denote by $\Delta_{\sigma}$ the stellar subdivision
of $\Delta$ on $\sigma$, by $x_{\sigma}$ the square-free monomial in $k[
\Delta]$ with support $\sigma$ and by $J_{\sigma}$ the annihilator of the
principal ideal of $k[ \Delta]$ generated by $x_{\sigma}$. Recall also from
\cite[Definition~1.2]{PR} that if $I= (f_{1},\dots,f_{r}) \subset R$ is a
homogeneous codimension $1$ ideal of a graded Gorenstein ring $R$ such that
the quotient $R/I$ is Gorenstein, then there exists $\phi\in\operatorname{Hom}%
_{R}(I,R)$ such that $\phi$ together with the inclusion $I \hookrightarrow R$
generate $\operatorname{Hom}_{R}(I,R)$ as an $R$-module. The Kustin--Miller
unprojection ring of the pair $I \subset R$ is defined as the quotient of
$R[y]$ by the ideal generated by the elements $yf_{i} - \phi(f_{i})$, where
$y$ is a new variable.

\begin{theorem}
\label{thm!mainalgebrathm} Suppose that $\Delta$ is a Gorenstein* simplicial
complex and that $\sigma\in\Delta$ is a face of dimension $d-1$ for some $d
\ge2$. Let $z$ be a new variable of degree $d-1$ and set $M =
\operatorname{Hom}_{k[ \Delta][z]} ((J_{\sigma}, z), k[ \Delta][z])$.

\begin{enumerate}
\itemsep=0pt

\item[\textrm{(a)}] $M$ is generated as a $k[ \Delta][z]$-module by the
elements $i$ and $\phi_{\sigma}$, where $i \colon(J_{\sigma}, z) \to k[
\Delta][z]$ is the natural inclusion morphism, and $\phi_{\sigma}$ is uniquely
specified by $\phi_{\sigma}(z) = x_{\sigma}$ and $\phi_{\sigma}(u) = 0 $ for
$u \in J_{\sigma}$.

\item[\textrm{(b)}] Denote by $S$ the Kustin--Miller unprojection ring of the
pair $(J_{\sigma}, z) \subset k[ \Delta][z]$. Then $z$ is a $S$-regular
element and $k[ \Delta_{\sigma}]$ is isomorphic to $S/(z)$ as a $k$-algebra.
\end{enumerate}
\end{theorem}

An example demonstrating Theorem~\ref{thm!mainalgebrathm} is the following.
Assume $\Delta$ is the boundary simplicial complex of the $2$-simplex and
$\sigma$ is a facet of $\Delta$. In coordinates, $k[\Delta] = k[x_{1}%
,x_{2},x_{3}]/(x_{1}x_{2}x_{3})$, $\sigma=\{1,2\}$ and $J_{\sigma} = 0
:(x_{1}x_{2}) = (x_{3})$. Then
\[
S = \frac{k[x_{1},\dots,x_{4},z]} {(x_{4}z-x_{1}x_{2}, x_{4}x_{3})},
\]
where $x_{4}$ denotes the new unprojection variable. Notice that when $z=0$,
$S|_{z=0}$ is isomorphic to $k[\Delta_{\sigma}]$, while when $a \in k^{*}$,
$S|_{z=a}$ is isomorphic (as ungraded $k$-algebra) to $k[\Delta]$. A toric
face ring interpretation of $S$ is discussed in
Example~\ref{example_with_picture}.

The paper is organised as follows: Theorem \ref{thm!mainalgebrathm} is proved
in Section \ref{sec!proof}. Section \ref{sec:pre} includes some definitions
and background related to the concepts which appear in Theorem
\ref{thm!mainalgebrathm}. Section \ref{sec!toric_face_ring_interpretion}
contains an interpretation of Theorem~\ref{thm!mainalgebrathm} using the
theory of toric face rings. In
Section~\ref{sec!application_to_stacked_polytopes}, we apply
Theorem~\ref{thm!mainalgebrathm} to inductively calculate the minimal graded
free resolutions of the Stanley--Reisner rings of the boundary simplicial
complexes of stacked polytopes, which were originally given in \cite{HM}. The
graded Betti numbers of these rings were first calculated in \cite{TH2}. When
the parameter value $d$ is not $3$, our methods allow us to obtain these Betti
numbers without using Hochster's formula or Alexander duality. We conclude in
Section \ref{sec:rem} with some remarks and directions for future research.

The applications of unprojection theory to Stanley--Reisner rings are not
limited to the case of stellar subdivisions, and in the paper \cite{BPcyclic}
we use unprojection techniques for an inductive treatment of Stanley--Reisner
rings associated to cyclic polytopes.

\section{Preliminaries}

\label{sec:pre}

Let $m$ be a positive integer and set $E = \{1, 2,\dots,m\}$. An (abstract)
\emph{simplicial complex} on the vertex set $E$ is a collection $\Delta$ of
subsets of $E$ such that (i) all singletons $\{i\}$ with $i \in E$ belong to
$\Delta$ and (ii) $\sigma\subset\tau\in\Delta$ implies $\sigma\in\Delta$. The
elements of $\Delta$ are called \emph{faces} and those maximal with respect to
inclusion are called \emph{facets}. The dimension of a face $\sigma$ is
defined as one less than the cardinality of $\sigma$. The dimension of
$\Delta$ is the maximum dimension of a face. The complex $\Delta$ is called
\emph{pure} if all facets of $\Delta$ have the same dimension. Any abstract
simplicial complex $\Delta$ has a geometric realization, which is unique up to
linear homeomorphism. When we refer to a topological property of $\Delta$, we
mean the corresponding property of the geometric realization of $\Delta$.

For any subset $\rho$ of $E$, we denote by $x_{\rho}$ the square-free monomial
in the polynomial ring $k[x_{1},\dots,x_{m}]$ with support $\rho$. The ideal
$I_{\Delta}$ of $k[x_{1},\dots,x_{m}]$ which is generated by the square-free
monomials $x_{\rho}$ with $\rho\notin\Delta$ is called the
\emph{Stanley-Reisner ideal} of $\Delta$. The \emph{face ring}, or
\emph{Stanley-Reisner ring}, $k[\Delta]$ of $\Delta$ over $k$, is defined as
the quotient ring of $k[x_{1},\dots,x_{m}]$ by the ideal $I_{\Delta}$. For a
face $\sigma$ of $\Delta$ denote by $\mathrm{lk}_{\Delta}(\sigma)=\{\tau
:\tau\cup\sigma\in\Delta$, $\tau\cap\sigma=\emptyset\}$ the \emph{link}, and
by $\mathrm{star}_{\Delta}(\sigma)=\{\tau:\tau\cup\sigma\in\Delta\}$ the
\emph{star} of $\sigma$ in $\Delta$. Given a face $\sigma$ of $\Delta$ of
dimension at least $1$, the \emph{stellar subdivision} of $\Delta$ on $\sigma$
is the simplicial complex $\Delta_{\sigma}$ on the vertex set $E\cup\{m+1\}$
obtained from $\Delta$ by removing all faces containing $\sigma$ and adding
all sets of the form $\tau\cup\{m+1\}$, where $\tau\in\Delta$ does not contain
$\sigma$ and $\tau\cup\sigma\in\Delta$. The complex $\Delta_{\sigma}$ is
homeomorphic to $\Delta$. We denote by $J{_{\sigma}}$ the ideal $(0:(x_{\sigma
}))$ of $k[\Delta]$, in other words
\[
J{_{\sigma}}=\{y\in k[\Delta]\colon yx_{\sigma}=0\}.
\]

The complex $\Delta$ is said to be Gorenstein* (over $k$) if $k[\Delta]$ is a
Gorenstein ring and given a vertex $i$ of $\Delta$ there 
exists $\sigma \in \Delta$ such that  $\sigma \cup \{ i \}$ is
not a face of $\Delta$.

 It is known
\cite[Section II.5]{Sta} that $\Delta$ is Gorenstein* if and only if for any
$\sigma\in\Delta$ (including the empty face) we have
\begin{equation}
\widetilde{H}_{i}(\mathrm{lk}_{\Delta}(\sigma),k)\ \cong\
\begin{cases}
k, & \text{if $i=\dim(\mathrm{lk}_{\Delta}(\sigma))$}\\
0, & \text{otherwise,}%
\end{cases}
\label{eq:gorenstein}%
\end{equation}
where $\widetilde{H}_{\ast}(\mathrm{lk}_{\Delta}(\sigma),k)$ denotes
simplicial homology of $\mathrm{lk}_{\Delta}(\sigma)$ with coefficients in the
field $k$. By \cite[Corollary~5.1.5]{BH}, any Gorenstein* complex $\Delta$ is
pure. It follows from (\ref{eq:gorenstein}) that the Gorenstein* property is
inherited by links. In particular, any codimension $1$ face of $\Delta$ is
contained in exactly $2$ facets of $\Delta$. The class of Gorenstein*
complexes includes all triangulations of spheres.

Assume $R$ is a polynomial ring over a field $k$ with the degrees of all
variables positive, and $M$ is a finitely generated graded $R$-module. Let
\[
0 \to F_{g} \to F_{g-1} \to\dots\to F_{1} \to F_{0} \to M\to0
\]
be the minimal graded free resolution of $M$ as $R$-module. Write
\[
F_{i} = \oplus_{j} R(-j)^{b_{ij}},
\]
then $b_{ij}$ is called the \emph{$ij$-th graded Betti number} of $M$, and we
also denote it by $b_{ij}(M)$. For more details about free resolutions and
Betti numbers see, for example, \cite[Sections~19, 20]{Ei}.

Assume $R$ is a ring. An element $r \in R$ will be called \emph{$R$-regular}
if the multiplication by $r$ map $R \to R, u \mapsto ru$ is injective. A
sequence $r_{1}, \dots,r_{n}$ of elements of $R$ will be called a
\emph{regular $R$-sequence} if $r_{1}$ is $R$-regular, and, for $2 \leq i \leq
n$, we have that $r_{i}$ is $R/(r_{1}, \dots,r_{i-1})$-regular.

\section{Proof of Theorem \ref{thm!mainalgebrathm}}

\label{sec!proof}

In this section, $\Delta$ denotes an $(n-1)$-dimensional simplicial complex on
the vertex set $\{1, 2,\dots,m\}$.

\begin{remark}
\label{rem!idealquotientismonomial} We will use the fact that $k[\Delta]$ has
no nonzero nilpotent elements and that if $I_{1},I_{2}$ are monomial ideals of
$k[\Delta]$, then so is the ideal quotient
\[
(I_{1}:I_{2})\ =\ \{y\in k[\Delta]\colon yI_{2}\subset I_{1}\}.
\]

\end{remark}

\begin{remark}
\label{rem!separatednessproperty} Assume that $\Delta$ is Gorenstein*. If $e$
is a vertex of $\Delta$ and $\sigma\in\Delta$ is a face that does not contain
$e$, then there exists a facet of $\Delta$ that contains $\sigma$ but not $e$.
Indeed, let $\tau_{1}$ be a facet of $\Delta$ containing $\sigma$. If
$\tau_{1}$ contains $e$, then there exists a facet $\tau_{2}$ distinct from
$\tau_{1}$ containing $\tau_{1}\setminus\{e\}$. This facet contains $\sigma$
and does not contain $e$.
\end{remark}

\begin{proposition}
\label{prop!technicalheart} Let $\Delta$ be a Gorenstein* simplicial complex
on the vertex set $\{1, 2,\dots,m\}$ and let $\sigma$ be a face of $\Delta$ of
dimension at least $1$. The ideal $J_{\sigma}$ is a codimension $0$ ideal of
$k[ \Delta]$ and the quotient $k[ \Delta]/J_{\sigma}$ is Gorenstein.
Moreover,
\[
(0 : J_{\sigma} ) = (x_{\sigma}).
\]

\end{proposition}

\begin{proof}
The first claim is well-known, cf. \cite[Theorem 21.23]{Ei}, and the second
follows from the observation that $k[\Delta]/J_{\sigma}=k[x_{1},\dots
,x_{m}]/I$, where
\[
I=I_{\mathrm{star}_{\Delta}(\sigma)}+(x_{i}:i\text{ is not a vertex of
}\mathrm{star}_{\Delta}(\sigma))\text{,}%
\]
and the fact that $\mathrm{lk}_{\Delta}(\sigma)$ is also Gorenstein*.

We now prove that $(0 : J_{\sigma}) = (x_{\sigma})$. It is clear that
$(x_{\sigma}) \subset(0 : J_{\sigma})$. Since $(0 : J_{\sigma})$ is a monomial
ideal (Remark \ref{rem!idealquotientismonomial}), it suffices to show that for
any nonzero monomial $u \in(0 : J_{\sigma})$ we have $u \in(x_{\sigma})$. Let
$\rho\in\Delta$ be the support of $u$. By the way of contradiction, suppose
that $u$ is not in $(x_{\sigma})$, so we may choose $i \in(\sigma\setminus
\rho)$. By Remark~\ref{rem!separatednessproperty}, there exists a facet $\tau$
of $\Delta$ which does not contain $i$ and contains $\rho$. Since $i$ is not
in $\tau$ and $\tau$ is a facet, we have $x_{i} x_{\tau}= 0$ in $k[ \Delta]$
and hence $x_{\tau}\in J_{\sigma}$. This fact and the assumption $u \in(0 :
J_{\sigma})$ imply that $x_{\tau}u = 0$ in $k[ \Delta]$. Since each variable
which appears in $u$ also appears in $x_{\tau}$, we conclude that $x_{\tau}$
is a nonzero nilpotent element of $k[ \Delta]$. This contradicts
Remark~\ref{rem!idealquotientismonomial} and completes the proof of the
proposition. \medskip
\end{proof}

\begin{remark}
The conclusion of Proposition~\ref{prop!technicalheart} is not true under the
weaker hypothesis that $k[\Delta]$ is Gorenstein. For a counterexample
consider
\[
\Delta=\{\{1,2\},\{1,3\},\{1\},\{2\},\{3\},\emptyset\}
\]
and $\sigma=\{1,2\}$. We have $k[\Delta]=k[x_{1},x_{2},x_{3}]/(x_{2}x_{3})$,
$J_{\sigma}=(0:x_{1}x_{2})=(x_{3})$, but $(0:J_{\sigma})=(x_{2})$. We believe
that this is also a counterexample to the second claim of Part a) of
\cite[Theorem 21.23]{Ei}, this is the reason we did not use this claim in the
proof of Proposition~\ref{prop!technicalheart}.
\end{remark}

Let $\sigma\in\Delta$ be a face of dimension $d-1$ for some $d \ge2$. We
recall that the stellar subdivision $\Delta_{\sigma}$ of $\Delta$ on $\sigma$
is a simplicial complex on the vertex set $\{1, 2,\dots,m+1\}$. We will use
the (easy) fact that
\begin{equation}
\label{eqn!aboutstellar}k[ \Delta_{\sigma}] \cong\frac{ k[x_{1},\dots,x_{m+1}]
} {(I_{ \Delta}, x_{\sigma}, x_{m+1} u_{1},\dots,x_{m+1}u_{r})},
\end{equation}
where $\{u_{1},\dots,u_{r}\}$ is a generating set of monomials for the ideal
$J_{\sigma}$ of $k[ \Delta]$.

\bigskip\noindent\emph{Proof of Theorem \ref{thm!mainalgebrathm}.} Clearly
there exists a unique element $\phi_{\sigma}$ of $M$ satisfying $\phi_{\sigma
}(z) = x_{\sigma}$ and $\phi_{\sigma}(u) = 0$ for $u \in J_{\sigma}$. Given $f
\in M$, we write $f(z) = w_{1} z + w_{2}$ with $w_{1} \in k[ \Delta][z]$ and
$w_{2} \in k[ \Delta]$ and set $g = f -w_{1} i \in M$, so that $g(z) = w_{2}$.
For $u \in J_{\sigma}$ we have
\[
z g(u) = g(zu) = u g(z) = uw_{2} \in k[ \Delta].
\]
Hence $g(u) = 0$ for all $u \in J_{\sigma}$, which implies $w_{2} \in(0 :
J_{\sigma})$. By Proposition~\ref{prop!technicalheart} we have $(0 :
J_{\sigma}) = (x_{\sigma})$. As a consequence, there exist $w \in k[ \Delta]$
such that $w_{2} = w x_{\sigma}$ and hence $g = w \phi_{\sigma}$. This proves
part (a) of the theorem.

By Proposition~\ref{prop!technicalheart}, the ring $k[ \Delta]/J_{\sigma}$ is
Gorenstein of the same dimension as $k[ \Delta]$. Therefore $(J_{\sigma},z)$
is a codimension $1$ homogeneous ideal of the graded Gorenstein ring $k[
\Delta][z]$, so the general theory of \cite{PR} applies. Using part (a) we
get
\[
S \cong\frac{ k[x_{1},\dots,x_{m+1},z] } {(I_{ \Delta}, x_{m+1}z-x_{\sigma},
x_{m+1} u_{1},\dots,x_{m+1} u_{r})},
\]
where the new variable $x_{m+1}$ has degree equal to $1$. It follows from
(\ref{eqn!aboutstellar}) that $S/(z) \cong k[ \Delta_{\sigma}]$. By
\cite[Theorem 1.5]{PR}, $S$ is Gorenstein of dimension equal to the dimension
of $k[ \Delta][z]$. As a consequence $\dim S/(z) = \dim S -1$ and therefore
$z$ is an $S$-regular element. This completes the proof of the theorem.
\qed \medskip

\section{Toric face ring interpretation}

\label{sec!toric_face_ring_interpretion}

Is it clear that Theorem~\ref{thm!mainalgebrathm} is equivalent to the
following theorem.

\begin{theorem}
\label{thm!mainthm_combinatorialformulation} Suppose that $\Delta$ is a
Gorenstein* simplicial complex and that $\sigma\in\Delta$ is a face of
dimension $d-1$ for some $d \ge2$. Let $z_{1}, \dots,z_{d-1}$ be $d-1$ new
variables of degree $1$ and set $M_{1} = \operatorname{Hom}_{k[ \Delta][z_{1},
\dots,z_{d-1}]} ((J_{\sigma}, z_{1}z_{2} \cdots z_{d-1}),$ $k[ \Delta][z_{1},
\dots,z_{d-1}])$.

\begin{enumerate}
\itemsep=0pt

\item[\textrm{(a)}] $M_{1}$ is generated as a $k[ \Delta][z_{1}, \dots
,z_{d-1}]$-module by the elements $i$ and $\phi_{\sigma}$, where $i
\colon(J_{\sigma}, z_{1}z_{2} \cdots z_{d-1}) \to k[ \Delta][z_{1},
\dots,z_{d-1}]$ is the natural inclusion morphism, and $\phi_{\sigma}$ is
uniquely specified by $\phi_{\sigma}(z_{1}z_{2} \cdots z_{d-1}) = x_{\sigma}$
and $\phi_{\sigma}(u) = 0 $ for $u \in J_{\sigma}$.

\item[\textrm{(b)}] Denote by $S_{1}$ the Kustin--Miller unprojection ring of
the pair \newline$(J_{\sigma}, z_{1}z_{2} \cdots z_{d-1}) \subset k[
\Delta][z_{1}, \dots,z_{d-1}]$. Then $z_{1},z_{2}, \dots,z_{d-1}$ is an
$S_{1}$-regular sequence, and $k[ \Delta_{\sigma}]$ is isomorphic to
$S_{1}/(z_{1},z_{2}, \dots,z_{d-1})$ as a $k$-algebra.
\end{enumerate}
\end{theorem}

We remark that, unlike in Theorem~\ref{thm!mainalgebrathm}, in
Theorem~\ref{thm!mainthm_combinatorialformulation} all variables have degree
$1$ which is the usual grading in the theory of Stanley--Reisner rings.
Compare also \cite[Section~4]{BPcyclic}, where a similar product $z_{1}z_{2}$
appears in a natural way when relating unprojection and cyclic polytopes.

Consider the Kustin--Miller unprojection ring
\[
S_{1} = \frac{ k[x_{1},\dots,x_{m+1},z_{1}, \dots,z_{d-1}] } {(I_{ \Delta},
x_{m+1}z_{1} \cdots z_{d-1}-x_{\sigma}, x_{m+1} u_{1},\dots,x_{m+1} u_{r})}
\]
appearing in Theorem~\ref{thm!mainthm_combinatorialformulation}, where as in
Section~\ref{sec!proof} $\{u_{1},\dots,u_{r}\}$ denotes a generating set of
monomials for the ideal $J_{\sigma} = (0 : x_{\sigma})$ of $k[ \Delta]$. We
will now give a combinatorial interpretation of $S_{1}$ using the notion of
toric face rings as defined by Stanley in \cite[p.~202]{Sta2}, compare also
\cite[Section 4]{BR} and \cite{BKR}. Let $M$ be a free $\mathbb{Z}$-module of
rank $m+d-1$, and consider the ${\mathbb{R}}$-vector space $M_{{\mathbb{R}}} =
M \otimes_{{\mathbb{Z}}} {\mathbb{R}}$. We will define a (finite, pointed)
rational polyhedral fan ${\mathcal{F}}$ in $M_{{\mathbb{R}}}$, such that
$S_{1}$ is isomorphic to the toric face ring $k[{\mathcal{F}}]$. For
simplicity of notation we assume in the following that $\sigma= \{1,2, \dots,d
\}$.

Denote by $e_{x,1}, \dots,e_{x,m}, e_{z,1}, \dots, e_{z,d-1} $ a fixed
${\mathbb{Z}}$-basis of $M$, and set
\[
e_{a} = (e_{x,1}+ \dots+ e_{x,d}) - (e_{z,1}+ \dots+ e_{z,d-1}) \in M.
\]
Assume $\tau=\{a_{1}, \dots,a_{p}\}$ is a face of $\Delta$. If $\sigma$ is not
a face of $\tau$ we set $c_{\tau}$ to be the cone in $M_{{\mathbb{R}}}$
spanned by the basis vectors
\[
e_{x,a_{1}}, \dots, e_{x,a_{p}}, e_{z,1}, \dots, e_{z,d-1},
\]
while if $\sigma$ is a face of $\tau$ we set $c_{\tau}$ to be the cone in
$M_{{\mathbb{R}}}$ spanned by the (non-affinely independent) vectors
\[
e_{x,a_{1}}, \dots, e_{x,a_{p}}, e_{z,1}, \dots, e_{z,d-1}, e_{a}.
\]
It is easy to see that the collection of cones $\{ c_{\tau} \bigm|  \tau\text{
face of } \Delta\}$ together with their faces form a fan ${\mathcal{F}}$ in
${M_{{\mathbb{R}}}}$ and that the toric face ring $k[{\mathcal{F}}]$ is
isomorphic as a $k$-algebra to $S_{1}$.

\begin{figure}[ptb]
\begin{center}
\includegraphics[
height=2.195in,
width=4.1597in,
trim= 0cm 16cm 0cm 3cm
]        {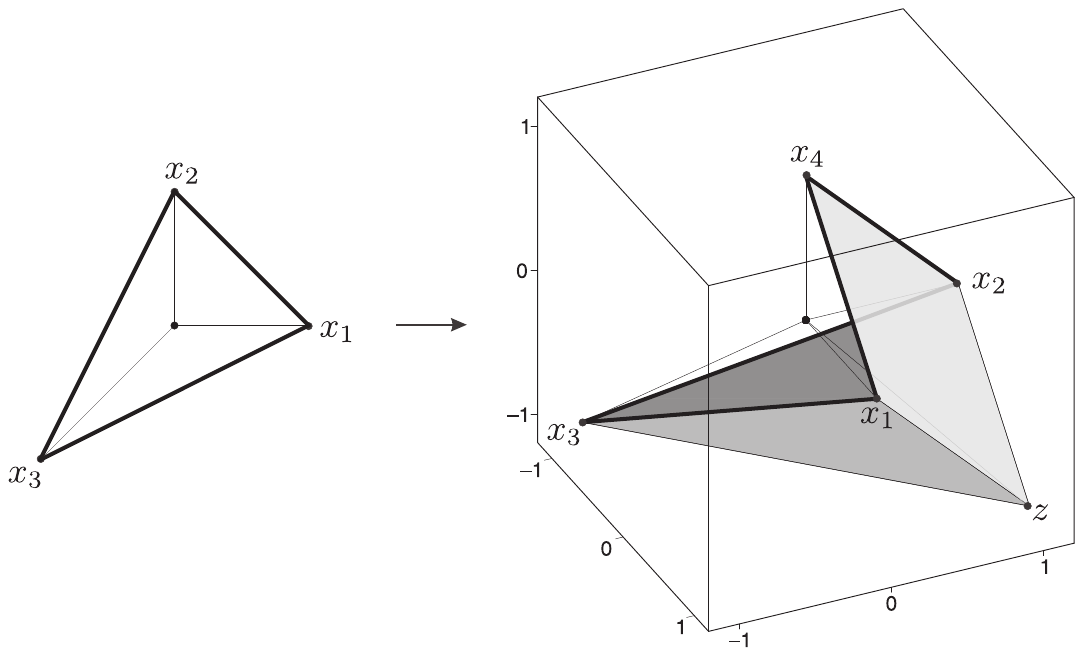}
\end{center}
\caption{Unprojection via toric face rings}%
\label{elliptic curve unprojection}%
\end{figure}

\begin{example}
\label{example_with_picture} Consider the example given after the statement of
Theorem~\ref{thm!mainalgebrathm}. That is, let $\Delta$ be the boundary of a
triangle with vertices corresponding to the variables $x_{1},x_{2},x_{3}$, and
denote by $\Delta_{\sigma}$ the stellar subdivision of $\Delta$ with respect
to the face $x_{1}x_{2}$. We embed the fan $\mathcal{F}$ into $\mathbb{R}^{3}$
by assigning to the variables $x_{1},x_{2},x_{3},x_{4}$ the rays generated by
$\left(  1,0,0\right)  $, $\left(  0,1,0\right)  $, $\left(  -1,-1,-1\right)
$, $\left(  0,0,1\right)  \in\mathbb{Z}^{3}$, i.e., those of the standard fan
of $\mathbb{P}^{3}$ as a toric variety. Then the ray associated to $z$ is
generated by $\left(  1,1,-1\right)  $. The right hand side of
Figure~\ref{elliptic curve unprojection} visualizes the Kustin-Miller
unprojection ring $S\cong k\left[  \mathcal{F}\right]  $ via representing each
cone of the embedded fan $\mathcal{F}$ by a polytope spanning it. There are
$3$ polytopes of maximal dimension, spanned by $\{x_{1},x_{3},z\},
\{x_{2},x_{3},z\}$ and $\{x_{1},x_{4},x_{2},z\}$. Notice that subdividing the
cone corresponding to $x_{1},x_{4},x_{2},z$ into $x_{1},x_{4},z$ and
$x_{4},x_{2},z$ amounts to passing from $S$ to the polynomial ring in the
variable $z$ over $k\left[  \Delta_{\sigma}\right]  $.
\end{example}

\section{Application to stacked polytopes}

\label{sec!application_to_stacked_polytopes}

\subsection{The Kustin--Miller complex construction}

\label{subs!generalities_about_KM_complexes}

The following construction, which is due to Kustin and Miller \cite{KM}, will
be important in the applications to stacked polytopes contained in
Subsection~\ref{subs!minimal_resol_for_stacked}.

Assume $R$ is a polynomial ring over a field with the degrees of all variables
positive, and $I\subset J\subset R$ are two homogeneous ideals of $R$ such
that both quotient rings $R/I$ and $R/J$ are Gorenstein and $\dim R/J=\dim
R/I-1$. We define $k_{1},k_{2}\in\mathbb{Z}$ such that $\omega_{R/I}%
=R/I(k_{1})$ and $\omega_{R/J}=R/J(k_{2})$, compare~\cite[Proposition~3.6.11]%
{BH}, and assume that $k_{1}>k_{2}$. Moreover, let
\[
0\rightarrow A_{g}\rightarrow A_{g-1}\rightarrow\dots\rightarrow
A_{1}\rightarrow A_{0}\rightarrow R/J\rightarrow0
\]
and
\[
0\rightarrow B_{g-1}\rightarrow\dots\rightarrow B_{1}\rightarrow
B_{0}\rightarrow R/I\rightarrow0
\]
be the minimal graded free resolutions of $R/J$ and $R/I$ respectively as
$R$-modules. Denote by $S=R[T]/Q$ the Kustin--Miller unprojection ring of the
pair $J\subset R/I$, where $T$ is a new variable of degree $k_{1}-k_{2}$.
Kustin and Miller constructed in \cite{KM} a graded free resolution of $S$ as
$R[T]$-module of the form
\[
0\rightarrow F_{g}\rightarrow F_{g-1}\rightarrow\dots\rightarrow
F_{1}\rightarrow F_{0}\rightarrow S\rightarrow0,
\]
where, when $g\geq3$,
\begin{align*}
&  F_{0}=B_{0}^{\prime},\quad\quad F_{1}=B_{1}^{\prime}\oplus A_{1}^{\prime
}(k_{2}-k_{1}),\\
\  &  F_{i}=B_{i}^{\prime}\oplus A_{i}^{\prime}(k_{2}-k_{1})\oplus
B_{i-1}^{\prime}(k_{2}-k_{1}),\quad\quad\text{for}\;\;2\leq i\leq g-2,\\
&  F_{g-1}=A_{g-1}^{\prime}(k_{2}-k_{1})\oplus B_{g-2}^{\prime}(k_{2}%
-k_{1}),\quad\quad F_{g}=B_{g-1}^{\prime}(k_{2}-k_{1}),
\end{align*}
cf. \cite[p.~307, Equation (3)]{KM}. When $g=2$ we have
\[
F_{0}=B_{0}^{\prime},\quad F_{1}=A_{1}^{\prime}(k_{2}-k_{1}),\quad F_{2}%
=B_{1}^{\prime}(k_{2}-k_{1}).
\]
In the above expressions, for an $R$-module $M$ we denoted by $M^{\prime}$ the
$R[T]$-module $M\otimes_{R}R[T]$. This resolution is, in general, not minimal,
see Example~\ref{example!nonminal_km_complex_construction} below. However, in
the case of stacked and cyclic polytopes it is minimal, see
Theorem~\ref{thm!applicationtostackedpolytopes} and \cite{BPcyclic}. We call
the complex consisting of the $F_{i}$ the \emph{Kustin--Miller complex
construction}. For more details and an implementation of this construction see
\cite{BP3}.

\begin{example}
Assume $k_{2}-k_{1}=-1$, and that the $2$ complexes are
\[
0\rightarrow A_{4}\rightarrow A_{3}\rightarrow A_{2}\rightarrow A_{1}%
\rightarrow A_{0}\rightarrow0
\]
and
\[
0\rightarrow B_{3}\rightarrow B_{2}\rightarrow B_{1}\rightarrow B_{0}%
\rightarrow0
\]
Then, the Kustin--Miller complex construction is of the form
\begin{multline*}
0\rightarrow B_{3}^{\prime}(-1)\rightarrow B_{2}^{\prime}(-1)\oplus
A_{3}^{\prime}(-1)\rightarrow B_{1}^{\prime}(-1)\oplus A_{2}^{\prime
}(-1)\oplus B_{2}^{\prime}\\
\rightarrow A_{1}^{\prime}(-1)\oplus B_{1}^{\prime}\rightarrow B_{0}^{\prime
}\rightarrow0
\end{multline*}
\textrm{ }
\end{example}

\begin{example}
\label{example!nonminal_km_complex_construction} Let $\Delta$ be the
simplicial complex with Stanley--Reisner ideal $(x_{1}x_{2}x_{3},x_{4}x_{5})$,
$\Delta$ is just the stellar subdivision of a facet of the boundary complex of
the $3$-simplex.  Then $\sigma=\{1,2\}$ is a face of $\Delta$. Since the Stanley-Reisner 
ideal  of $\Delta_{\sigma}$  is minimally generated by $3$  monomials and 
not by $5$, the Kustin--Miller complex construction  gives a graded resolution 
of $k[\Delta_{\sigma}]$ which is not minimal.
\end{example}

\subsection{ The minimal resolution for stacked polytopes}

\label{subs!minimal_resol_for_stacked}

Assume $d\geq2$ is a fixed integer. Recall from \cite[p.~448]{TH2}, that
starting from a $d$-simplex one can add new vertices by building shallow
pyramids over facets to obtain a simplicial convex $d$-polytope with $m$
vertices, called a \emph{stacked polytope} $P_{d}(m)$. We denote by $\Delta
P_{d}(m)$ the boundary simplicial complex of the simplicial polytope
$P_{d}(m)$. By definition, $\Delta P_{d}(m)$ has as elements the empty set and
the sets of vertices of the proper faces of $P_{d}(m)$,
cf.~\cite[Corollary~5.2.7]{BH}. There is a slight abuse of notation here,
since the combinatorial type of $\Delta P_{d}(m)$ does not depend only on $d$
and $m$ but also on the specific choices of the sequence of facets we used
when building the shallow pyramids. The graded Betti numbers $b_{ij}$ of the
Stanley-Reisner ring $k[\Delta P_{d}(m)]$ have been calculated  by Terai and Hibi 
in \cite[Theorem~1.1]{TH2}, and it turns out that they only depend on $d$ and
$m$. Later Herzog and Li Marzi \cite{HM}  constructed the minimal graded free resolution
of $k[\Delta P_{d}(m)]$.  In Theorem~\ref{thm!applicationtostackedpolytopes}
we give a different proof of their result based on Theorem~\ref{thm!mainalgebrathm}.

It is clear that, for $d<m$, the simplicial complex $\Delta P_{d}(m+1)$ can be
considered as the stellar subdivision of the boundary simplicial complex
$\Delta P_{d}(m)$ of a stacked polytope $P_{d}(m)$ with respect to a facet
$\sigma$ of $\Delta P_{d}(m)$. Since $\sigma$ is a facet, the ideal
$(J_{\sigma},z)$ is generated by the regular sequence $x_{\rho},z$, where
$\rho$ takes values in the set of vertices of $\Delta P_{d}(m)$ which are not
vertices of $\sigma$. Hence, the minimal graded free resolution of
$(J_{\sigma},z)$ is a Koszul complex. Combining
Theorem~\ref{thm!mainalgebrathm} with the Kustin--Miller complex construction
described in Subsection~\ref{subs!generalities_about_KM_complexes} we can get,
starting with the Koszul complex and the minimal graded free resolution of
$k[\Delta P_{d}(m)]$, a graded free resolution of $k[\Delta P_{d}(m+1)]$. The
following theorem states that we indeed get the minimal graded free resolution
of $k[\Delta P_{d}(m+1)]$. In this way we recover the result from \cite{HM}
using different ideas. We remark that, when $d=2$ or $d\geq4$, we do not use
in the proof of the theorem the calculation of the graded Betti numbers of
$k[\Delta P_{d}(m)]$ given in \cite{TH2}, and, moreover, we obtain these
numbers in Proposition~\ref{prop!keytechnicalfor_stacked_application}. The
proof of the theorem will be given in
Subsection~\ref{subs!proofofapplictostacked}.

\begin{theorem}
\label{thm!applicationtostackedpolytopes} Assume $d\geq2$ and $d+1<m$. The
resolution of $k[\Delta P_{d}(m+1)]$, obtained using the Kustin--Miller
complex construction starting from the minimal graded free resolution of
$k[\Delta P_{d}(m)]$ and the Koszul complex resolving $(J_{\sigma},z)$ is minimal.
\end{theorem}

\subsection{Proof of Theorem~\ref{thm!applicationtostackedpolytopes}}

\label{subs!proofofapplictostacked}

We need the following combinatorial definition. Assume $d \geq2$ and $d < m$.
For $1 \leq i \leq m-d-1$ we define
\[
\theta(d,m,i) = i \binom{m-d}{i+1},
\]
compare~\cite[p.~448]{TH2}. Moreover we set $\theta(d,m,0) = \theta(d,m,m-d) =
0$.


\begin{lemma}
\label{lem!combinatorial_induction_for_stacked} \textrm{(Compare
\cite[p.~451]{TH2}).} Assume $1 \leq i \leq m-d$. Then
\[
\theta(d,m+1,i) = \theta(d,m,i) + \binom{m-d}{i} +\theta(d,m,i-1).
\]
(By our conventions, for $i=1$ the equality becomes $\theta(d,m+1,1) =
\theta(d,m,1) + (m-d)$, while for $i = m-d$ it becomes $\theta(d,m+1,m-d) =
\theta(d,m,m-d-1) + 1$).
\end{lemma}

\begin{proof}
Assume first $2 \leq i \leq m-d-1$. Using the Pascal triangle identity
$\binom{m}{d} = \binom{m-1}{d} + \binom{m-1}{d-1}$ we have {\small
\begin{align*}
\theta(d,m+1,i)  &  = i \binom{m+1-d}{i+1} = i (\binom{m-d}{i+1} + \binom
{m-d}{i})\\
&  = i \binom{m-d}{i+1} + \binom{m-d}{i} + (i-1) \binom{m-d}{i}\\
&  = \theta(d,m,i) + \binom{m-d}{i} + \theta(d,m,i-1).
\end{align*}
} The special cases $i=1$ and $i=m-d$ are proven by the same argument.
\end{proof}

The following proposition is well-known.

\begin{prop}
\label{prop!dividing_resolutions_by_graded_elements}
\textrm{(\cite[Proposition~1.1.5]{BH}).} Assume $R= k[x_{1}, \dots,x_{n}]$ is
a polynomial ring over a field $k$ with the degrees of all variables positive,
and $I \subset R$ a homogeneous ideal. Moreover, assume that $x_{n}$ is
$R/I$-regular. Denote by $cF$ the minimal graded free resolution of $R/I$ as
$R$-module. We then have that $cF \otimes_{R} R/(x_{n})$ is the minimal graded
free resolution of $R/(I,x_{n})$ as $k[x_{1}, \dots,x_{n-1}]$-module, where we
used the natural isomorphisms \newline$R \otimes_{R} R/(x_{n}) \cong R/(x_{n})
\cong k[x_{1}, \dots,x_{n-1}]$.
\end{prop}

The proof of the following proposition is an immediate corollary of the
construction of the Koszul complex in \cite[Section~1.6]{BH}.

\begin{prop}
\label{prop!Koszul_complex_with_different_degrees} Assume $R= k[x_{1},
\dots,x_{n}]$ is a polynomial ring over a field $k$ with the degrees of all
variables positive, $p \leq n$ a fixed integer, and $g_{1}, \dots,g_{p}$, an
$R$-regular sequence consisting of homogeneous elements of $R$, with $\deg
g_{i} = 1$, for $1 \leq i \leq p-1$, and $\deg g_{p} = q \geq1$. Then, the
minimal resolution of $R/(g_{1}, \dots,g_{p})$ is of the form
\[
0 \to F_{p} \to F_{p-1} \dots\to F_{1} \to F_{0},
\]
with $F_{0} \cong R$, $F_{p} \cong R(-p-q+1)$, and
\[
F_{i} \cong R(-i)^{b_{i}} \oplus R(-q-i+1)^{b_{p-i}}
\]
for $1 \leq i \leq p-1$, where $b_{i} = \binom{p-1}{i}$.
\end{prop}

Theorem~\ref{thm!applicationtostackedpolytopes} is an immediate consequence of 
the following
more precise proposition. Notice that as we have already mentioned the
statements about the graded Betti numbers of $k[\Delta P_{d}(m)]$ have been
proven before in \cite{TH2}. For $d \not = 3$ we do not use in our proof the
results of \cite{TH2}.

\begin{prop}
\label{prop!keytechnicalfor_stacked_application} Assume $d \geq2$ and $d+1 <
m$. Set $b_{ij} = b_{ij}(k[\Delta P_{d}(m)])$. Then the statement of
Theorem~\ref{thm!applicationtostackedpolytopes} is true for $(d,m)$. Moreover,
we have that if $d=2$ then $b_{ij}= 1$ for $(i,j) \in\{ (0,0), (m-d,m) \}$,
\[
b_{i,i+1} = \theta(d,m,i) + \theta(d,m,m-d-i),
\]
for $1 \leq i \leq m-d-1$, and $b_{ij} = 0$ otherwise. If $d \geq3$, we have
$b_{ij}= 1$ for $(i,j) \in\{ (0,0), (m-d,m) \}$,
\[
b_{i,i+1} = \theta(d,m,i), \quad b_{i,d+i-1} = \theta(d,m,m-d-i),
\]
for $1 \leq i \leq m-d-1$, and $b_{ij} = 0$ otherwise.
\end{prop}

\begin{proof}
We fix $d \geq2$ and use induction on $m$. If $d \geq2$ and $m = d+2$ then
$k[\Delta P_{d}(m)]$ is a type $(2,d)$ codimension $2$ complete intersection
and everything is clear.

Assume $d \not = 3$, and that
Proposition~\ref{prop!keytechnicalfor_stacked_application} is true for
$(d,m)$. By Theorem~\ref{thm!mainalgebrathm}, the extension ring $S$ of
$k[\Delta P_{d}(m+1)]$ is the Kustin--Miller unprojection ring of the pair
$(J_{\sigma},z) \subset k[\Delta P_{d}(m)][z]$. As we noticed above, the ideal
$(J_{\sigma},z)$ is generated by a regular sequence, so the Koszul complex
described in Proposition~\ref{prop!Koszul_complex_with_different_degrees} is
the minimal resolution of $k[\Delta P_{d}(m)][z]/(J_{\sigma}, z)$. Combining
Proposition~\ref{prop!dividing_resolutions_by_graded_elements} and the
discussion of Subsection~\ref{subs!generalities_about_KM_complexes}, starting
from the Koszul complex and the minimal graded free resolution of $k[\Delta
P_{d}(m)]$, the Kustin--Miller complex construction gives a graded free
resolution of $k[\Delta P_{d}(m+1)]$. Using
Lemma~\ref{lem!combinatorial_induction_for_stacked} this complex has the
conjectured graded Betti numbers, and since there are no degree $0$ morphisms
it is necessarily minimal.

When $d=3$ the above arguments work except for the minimality argument, since
there are degree $0$ morphisms. But comparing the graded Betti number of the
Kustin--Miller complex construction with the graded Betti numbers of $k[\Delta
P_{d}(m+1)]$ calculated in \cite{TH2} we again obtain the minimality of the
Kustin--Miller complex construction.
\end{proof}

\section{Remarks and open questions}

\label{sec:rem}

In the following we use the notation introduced in Section \ref{intro} and
Theorem \ref{thm!mainalgebrathm}.

\begin{remark}
It follows from Theorem~\ref{thm!mainalgebrathm} that $S$ is a $1$-parameter
deformation ring of $k[\Delta_{\sigma}]$, compare \cite[Exerc.~18.18]{Ei}. The
fact that such a deformation ring of $k[\Delta_{\sigma}]$ exists is a special
case of more general results due to Altmann and Christophersen \cite{AC1,AC2}.
\end{remark}

\begin{remark}
Using the Kustin--Miller complex construction described in
Subsection~\ref{subs!generalities_about_KM_complexes}, we can construct a
graded free resolution of $S$, therefore using
Proposition~\ref{prop!dividing_resolutions_by_graded_elements} also of
$k[\Delta_{\sigma}]$, starting from graded free resolutions of $k[\Delta]$ and
$k[\Delta]/J_{\sigma}$. In particular, it follows that
\[
F(k[\Delta_{\sigma}],t)=F(k[\Delta],t)+(t+t^{2}+\cdots+t^{d-1})\,F(k[\Delta
]/J_{\sigma},t),
\]
where $F(R,t)$ stands for the Hilbert series of $R$ and $d-1$ is the dimension
of the face $\sigma$. This equality can be rewritten as
\begin{equation}
h(\Delta_{\sigma},t)=h(\Delta,t)+(t+t^{2}+\cdots+t^{d-1})\,h(lk_{\Delta
}(\sigma),t), \label{eq:h}%
\end{equation}
where $h(\Gamma,t)$ stands for the $h$-polynomial \cite[Section II.2]{Sta} of
the simplicial complex $\Gamma$. It is not hard to see that (\ref{eq:h}) holds
for any pure simplicial complex $\Delta$. Indeed, one can check directly that
(\ref{eq:h}) is equivalent to the formula
\[
f_{j}(\Delta_{\sigma})=f_{j}(\Delta)-f_{j-d}(lk_{\Delta}(\sigma))\,+\,\sum
_{i\geq0}\ \binom{d}{j-1}f_{i-1}(lk_{\Delta}(\sigma)),
\]
where $f_{j}(\Gamma)$ denotes the number of $j$-dimensional faces of a complex
$\Gamma$. That formula follows from the definition of $\Delta_{\sigma}$.
\end{remark}

\begin{remark}
In \cite{NP}, Neves and the second author introduced the $\binom{n}{2}$
Pfaffians format, starting from a certain hypersurface ideal. We give a
monomial interpretation of the construction. Start with the boundary
simplicial complex $\Delta$ of the $(n-1)$-simplex. Denote by $\Delta_{1}$ the
simplicial complex obtained by the stellar subdivisions of all facets of
$\Delta$. It is easy to check that the Stanley--Reisner ideal of $\Delta_{1}$
is equal to $\widetilde{I_{n}}$, where $\widetilde{I_{n}}$ denotes the ideal
obtained by substituting $z_{i}=0$, for $1\leq i\leq n$, and $r_{d_{1}%
,\dots,d_{n}}=1$, for $(d_{1},\dots,d_{n})\in\{0,1\}^{n}$, to the ideal
$I_{n}$ defined in \cite[Definition~2.2]{NP}.

Similarly, in \cite[Section~4.3]{NP2}, Neves and the second author constructed
a codimension $11$ Gorenstein ideal starting from a certain codimension $2$
complete intersection ideal. The monomial interpretation of the construction
is as follows. Denote by $\Delta$ the simplicial complex which is the join
\cite[p.~221]{BH} of $2$ copies of the boundary simplicial complex of the
$2$-simplex. $\Delta$ has Stanley--Reisner ideal equal to $(x_{11}x_{12}%
x_{13},x_{21}x_{22}x_{23})$ and exactly $9$ facets. Denote by $\Delta_{1}$ the
simplicial complex obtained by the stellar subdivisions of $\Delta$ on these
$9$ facets. Using the notations of \cite[Section~2]{NP2}, denote by
$I_{\mathcal{L}}$ the kernel of the surjection $R[y_{u}\bigm|u\in
L]\rightarrow R_{\mathcal{L}}$. It is easy to check that the Stanley--Reisner
ideal of $\Delta_{1}$ is equal to $\widetilde{I_{\mathcal{L}}}$, where
$\widetilde{I_{\mathcal{L}}}$ denotes the ideal obtained by substituting
$x_{3i}=0$, for $1\leq i\leq3$, to $I_{\mathcal{L}}$.
\end{remark}

\begin{remark}
It is plausible that our ideas also generalize to non-Gorenstein simplicial
complexes. To do this a more detailed study of non-Gorenstein unprojections
would be necessary.
\end{remark}

\begin{remark}
Combining our results with those of \cite{KM3} we get a link between stellar
subdivisions of  Gorenstein* simplicial complexes and linkage theory
\cite{Mi}. Is it possible to use this connection to define new combinatorial
invariants of simplicial complexes?
\end{remark}

\vspace{0.1in} \noindent\emph{Acknowledgements}. The authors are grateful to
Christos Athanasiadis for important discussions and suggestions. They also
thank Jo\~{a}o Martins for useful discussions, Tim R\"{o}mer for useful
comments on an earlier version, and an anonymous referee
for informing us about \cite{HM}.

\end{document}